\documentclass[12pt]{amsart}
\usepackage{amsmath,amsthm,latexsym,amscd,amsbsy,amssymb}
 \relax

\chardef\bslash=`\\ 

\makeatletter
\def\verbatim{\interlinepenalty\@M \@verbatim
  \leftskip\@totalleftmargin\advance\leftskip2pc
  \frenchspacing\@vobeyspaces \@xverbatim}
\makeatother
\makeatletter
  \def\dgt@k{\dg@DX=-3 \dg@DY=2 \dg@SIZE=3}
\makeatother

\makeatletter
  \def\dgt@kk{\dg@DX=3 \dg@DY=-1 \dg@SIZE=3}%
\makeatother

\theoremstyle{plain}

\newtheorem*{A}{Theorem}
\newtheorem*{B}{Lemma}

\theoremstyle{definition}

\numberwithin{equation}{section}

\def\Int{{\rm Int}}

\def\ed{{\rm e{\text -}dim}}
\def\wh{\widehat}
\def\wt{\widetilde}
\def\St{{\rm St}}

\begin{document}

\title[Extension dimensional approximation theorem]
{Extension dimensional approximation theorem}
\author{N. Brodsky}
\address{Department of Mathematics and Statistics,
University of Saskatche\-wan,
McLean Hall, 106 Wiggins Road, Saskatoon, SK, S7N 5E6,
Canada}
\email{brodsky@math.usask.ca}
\author{Alex Chigogidze}
\address{Department of Mathematics and Statistics,
University of Saskatche\-wan,
McLean Hall, 106 Wiggins Road, Saskatoon, SK, S7N 5E6,
Canada}
\email{chigogid@math.usask.ca}
\thanks{Second author was partially supported by NSERC research grant.}


\begin{abstract}{It is known that if an upper semicontinuous multivalued mapping $F \colon X \to Y$, defined on an $n$-dimensional compactum $X$, has $UV^{n-1}$-point images, then every neighbourhood of the graph of $F$ (in the product $X \times Y$) contains the graph of a single-valued continuous mapping $f \colon X \to Y$. Similar result is known to be true when $X$ is a compact $C$-space and images of $F$ have trivial shape. We extend and unify both of these results in terms of extension theory.}
\end{abstract}

\maketitle
\markboth{N.~Brodsky, A.~Chigogidze}{Extension dimensional
approximation theorem}

\section{Introduction}\label{S:introdution}

Single-valued approximations of multivalued maps are proved to be
very useful in geometric topology, fixed point theory,
control theory and others (see a survey~\cite{Kr1}).
We consider the problem of single-valued continuous
graph-approximation of upper semicontinuous (u.s.c.)
multivalued mappings. We say that a multivalued mapping $F\colon X\to Y$
admits graph-approximations if every neighborhood
of the graph of $F$
(in the product $X \times Y$) contains the graph of a single-valued
continuous mapping $f\colon X\to Y$.

Essentially there are three types of results concerning our problem.
First assumes that multivalued mappings $F\colon X\to Y$
have $UV^{n-1}$ point-images and $\dim X\le n$
(see~\cite{L, ShB, Kr2}).
The second type of results deal with $UV^\infty$-valued
mappings defined on C-spaces~\cite{Anc}.
Finally results of the third type consider $UV^\infty$-valued
mappings defined on $ANR$-spaces~\cite{GGK, Kr2}.

In this paper we prove an approximation theorem which generalizes and
unifies the known results of the first and second types. Unification is
achieved by exploiting recently created \cite{Dr1}, \cite{Dr2}
theory of extension dimension
and associated to it concepts of homotopy and shape \cite{Ch}.
Precise definitions will be given below in Section \ref{S:preliminaries}.
Here we only provide some of the notation related to the extension
dimension.

Let $L$ be a CW-complex.
A space $X$ is said to have {\it extension dimension} $\le [L]$
(notation: $\ed X \le [L]$) if
any mapping of its closed subspace
$A\subset X$ into $L$ admits an extension to the whole space
$X$\footnote{Everywhere below $[L]$ denotes the class of
complexes generated by $L$ with respect to the above extension
property, see \cite{Dr1}, \cite{Dr2}, \cite{Ch} for details.}. It
is known that $\dim X \leq n$ is equivalent to
$\ed X \leq \left[ S^{n}\right]$ and that $\dim_{G}X \leq n$
is equivalent to $\ed X\leq \left[ K(G,n)\right]$ ($K(G,n)$ stands
for the corresponding Eilenberg-MacLane complex). One
can develop homotopy and shape theories
specifically designed to work for at most
$[L]$-dimensional spaces. Compacta of trivial
$[L]$-shape are precisely $UV^{[L]}$-compacta \cite{Ch}.

Now we are ready to formulate our main result.

\begin{A}\label{T:A}
Let $L$ be a countable CW-complex and
$F\colon X\to Y$ be an u.s.c. $UV^{[L]}$-valued
mapping of a paracompact space $X$ to a completely
metrizable space $Y$.
If $X$ is C-space of extension dimension $\ed X\le [L]$,
then every neighborhood of the graph of $F$ contains the graph
of a single-valued continuous mapping $f\colon X\to Y$.
\end{A}

Note that if $L$ is the sphere $S^n$,
we obtain an approximation theorem for
$UV^{n-1}$-valued mappings of $n$-dimensional space.
And if $L$ is a point (or any other contractible complex),
we obtain a theorem of Ancel
on approximations of $UV^\infty$-valued mappings
of C-space~\cite{Anc}.

What do we need to construct a mapping from a space $X$?
Suppose that we can construct and, moreover, extend a
mapping from $X$ locally. Then one can try to obtain
a fine cover of $X$ and to construct a global
mapping by induction, extending it successively over
"skeleta" of this cover. The problem is to control this
process when the cover has infinite order.
Property C gives us a possibility of such a control.

Let us explain this with a bit more detail.
A topological space $X$ has {\it property} C if for each sequence
$\{u^i\mid i\ge 1\}$ of open covers of $X$, there is an open cover
$\Sigma$ of $X$ of the form $\cup_{i=1}^\infty \sigma_i$
such that for each $i\ge 1$, $\sigma_i$ is a
pairwise disjoint collection which refines $u^i$.
If the space $X$ is paracompact, we can choose the cover
$\Sigma$ to be locally finite.
The cover $\Sigma$ has very important property that every
"simplex" $\{s_0,\dots,s_n\}$ of this cover (i.e. the set
of elements $\{s_0,\dots,s_n\}$ such that $s_0\cap\dots\cap s_n\ne\emptyset$)
has a natural order on its vertices. Indeed, for any
element $s\in\Sigma$ denote by $\sigma(s)$ the integer
such that $s\in \sigma_{\sigma(s)}$.
Since $s_i\cap s_j\ne\emptyset$, then $\sigma(s_i)\ne\sigma(s_j)$
and we can order elements $s_0, s_1,\dots,s_n$ according to
the order of numbers $\sigma(s_0),\sigma(s_1),\dots,\sigma(s_n)$.

We take a cover $\Sigma$ of $X$ which refines our fine
cover so that every simplex $\langle \sigma_0,\dots,\sigma_n\rangle $
of the nerve $N(\Sigma)$ has a natural order on its vertices.
Then every simplex has a {\it basic} vertex (merely the smallest one).
For every vertex $\langle \sigma\rangle $ of $N(\Sigma)$ (i.e. for every element
$\sigma$ of the cover $\Sigma$) we fix a "rule" of extension
of mappings defined on subsets of $\sigma$.
Then the process goes on by induction on dimension of "skeleta"
as follows: for a "simplex" the extension of mapping from
the "boundary" to the "interior" is induced by the rule
of the basic vertex of this simplex.
Obviously the mapping on a simplex depends only on the basic vertex
of this simplex and does not depend on the dimension of the simplex.
This provides the needed control.

\section{Preliminaries}\label{S:preliminaries}

Let us recall some definitions and introduce our notations.
We denote by $\Int A$ and $\overline A$ the interior and closure
of the set $A$ respectively.
For a cover $\omega$ of a space $X$ and for a subset $A \subseteq X$ let
$\St(A,\omega)$ denote the star of the set $A$ with respect to $\omega$.

The {\it graph} of a multivalued mapping $F\colon X\to Y$
is the subset $\Gamma_F=\{(x,y)\in X\times Y\colon y\in F(x)\}$
of the product $X\times Y$.
A multivalued mapping $F\colon X\to Y$ is called
{\it upper semicontinuous} (notation: u.s.c.) if for any open set
$U\subset Y$ the set $\{x\in X\colon F(x)\subset U\}$ is open in $X$.

Let $L$ be a CW-complex. A pair of spaces $V\subset U$ is said
to be {\it $[L]$-connected}
if for every paracompact space $X$ of extension dimension $\ed X\le [L]$
and for every closed subspace $A\subset X$ any mapping of $A$
into $V$ can be extended to a mapping of $X$ into $U$.
A compact subspace $K\subset Z$ is called {\it $UV^{[L]}$-compactum in $Z$}
if any neighborhood $U$ of $K$ contains a neighborhood $V$ of $K$
such that the pair $V\subset U$ is $L$-connected.
A compact-valued mapping $F\colon X\to Y$ is called
{\it $UV^{[L]}$-valued} if for any point $x\in X$
the set $F(x)$ is $UV^{[L]}$-compactum in $Y$.
A mapping $f\colon Y\to X$ is said to be {\it $[L]$-soft}
if for any paracompact space $Z$ with $\ed Z \le [L]$, its closed subspace
$A\subset Z$ and any mappings $g\colon Z\to X$ and
$\wt g_A\colon A\to Y$ such that $f\circ \wt g_A=g|_A$
there exists a mapping $\wt g\colon Z\to Y$ such that $f\circ \wt g=g$.
Finally let $AE([L])$ denote the class of spaces with
$[L]$-soft constant mappings.

Now we introduce the notion of $[L]$-extension which will represent
a "rule" for extending mappings in the proof of our theorem.
Let $V\subset U$ be a pair of spaces. An {\it $[L]$-extension}
of the space $V$ with respect to $U$ is a pair $V'\subset W$ of spaces
and a mapping $e\colon W\to U$ such that:

\begin{itemize}
\item[1)] $W\in AE([L])$;
\item[2)] $e|_{V'}$ is $[L]$-soft mapping onto $V$.
\end{itemize}

The following is a key property of $[L]$-extensions needed
in the proof (Section \ref{S:proof}) of our theorem.
Let a pair $V'\subset W$ of spaces and a mapping $e\colon W\to U$
represent an $[L]$-extension of the pair $V\subset U$.
\medskip

{\bf $[L]$-extension property.} {\em Let $A\subset B$ be a pair of closed
subspaces of paracompact space $X$ of extension dimension $\ed X \le [L]$.
Suppose that we have mappings $f\colon B\to U$ and $g\colon A\to W$
such that $e\circ g=f|_A$, $f(\overline{B\setminus A}) \subset V$ and $g(A\cap \overline{B\setminus A}) \subset V'$.
Then there exists a mapping $g'\colon X\to W$ such that $e\circ g'|_B=f$.}
\medskip

We construct $g'$ in two steps. First, we use $[L]$-softness
of $e$ over $V$ to extend $g$ to a mapping $\wt g\colon B\to W$
such that $e\circ \wt g=f$ (we apply $[L]$-softness to the $[L]$-dimensional
pair $A\cap \overline{B\setminus A}\subset \overline{B\setminus A}$).
Finally we can extend $\wt g$ to the space $X$ since $W$ is $AE([L])$.

\begin{B}\label{L:B}
Let $V\subset U$ be $[L]$-connected pair.
If $V$ is completely metrizable space,
then $V$ admits an $[L]$-extension with respect to $U$.
\end{B}

\begin{proof}
There exists a completely metrizable space $V'$ with $\ed V' \leq [L]$
and an $[L]$-soft mapping $e_V\colon V'\to V$~\cite{ChV}.
Consider an $AE([L])$-space $W$ of dimension $\ed W \leq [L]$
containing $V'$ as a closed subspace~\cite{ChV}.
Since the pair $V\subset U$ is $[L]$-connected, we can extend
the mapping $e_V$ to a mapping $e\colon W\to U$.
\end{proof}

\section{Proof of the theorem}\label{S:proof}

For a given $UV^{[L]}$-valued mapping $F\colon X\to Y$ we fix an arbitrary
neighborhood ${\mathcal U}\subset X\times Y$ of its graph $\Gamma_F$.
The proof of our theorem consists of the following two steps.
\bigskip

{\bf 1. Construction of families of rectangles.}
\medskip

For every integer $i\ge 0$ we construct families of open rectangles
$\{u^i_\lambda\times U^i_\lambda\}_{\lambda\in\Lambda_i}$ and
closed rectangles $\{v^i_\mu\times V^i_\mu\}_{\mu\in M_i}$
in the product $X\times Y$ such that:
\smallskip

\begin{itemize}
\item[(1)] $u^0_\lambda\times U^0_\lambda\subset {\mathcal U}$ for every
   $\lambda\in\Lambda_0$;
\smallskip

\item[(2)] $u^i=\{u^i_\lambda\}_{\lambda\in\Lambda_i}$ and
   $v^i=\{v^i_\mu\}_{\mu\in M_i}$ are coverings of $X$
   (in fact, $\{\Int\ v^i_\mu\}_{\mu\in M_i}$ are coverings of $X$);
\smallskip

\item[(3)] $F(u^i_\lambda)\subset U^i_\lambda$ and
   $F(v^i_\mu)\subset \Int\ V^i_\mu$
   for every $i\ge 0$, $\mu\in M_i$ and $\lambda\in\Lambda_i$;
\smallskip

\item[(4)] for every $i\ge 0$ and every $\mu\in M_i$ there exists $\lambda\in\Lambda_i$
   such that $V^i_\mu\subset U^i_\lambda$, $v^i_\mu\subset u^i_\lambda$,
   and the pair $V^i_\mu\subset U^i_\lambda$ is $[L]$-connected;
\smallskip

{\it Choice 1.} For given $i\ge 0$ and $\mu\in M_i$ we fix such a
$\lambda=\lambda(\mu)$, and for $[L]$-connected pair
$V^i_\mu\subset U^i_\lambda$, by Lemma, we can fix $[L]$-extension
$e^i_\mu\colon (\wt V^i_\mu, W^i_\mu)\to (V^i_\mu,U^i_\lambda)$;
\smallskip

\item[(5)] for every $i\ge 0$ and every $\lambda\in\Lambda_{i+1}$ there exists
   $\mu\in M_i$ such that $\St(u^{i+1}_\lambda,u^{i+1})\subset v^i_\mu$
   and every rectangle $u^{i+1}_\gamma\times U^{i+1}_\gamma$
   is contained in the rectangle $v^i_\mu\times V^i_\mu$
   provided $u^{i+1}_\gamma\cap u^{i+1}_\lambda\ne\emptyset$;
\smallskip

{\it Choice 2.} For given $i\ge 0$ and $\lambda\in\Lambda_{i+1}$
we fix such a $\mu=\mu(\lambda)$.
\end{itemize}

First, we construct a family
$\{u^0_\lambda\times U^0_\lambda\}_{\lambda\in\Lambda_0}$.
Put $\Lambda_0= X$ and for a point $x\in X$ consider a rectangle
$u_x\times U_x^0\subset {\mathcal U}$ such that $F(x)\subset U_x^0$
(existence of such a rectangle follows from compactness of $F(x)$).
Since $F$ is u.s.c., we can choose a neighborhood $u_x^0\subset u_x$
of the point $x$ such that $F(u_x^0)\subset U^0_x$.

The construction of families of rectangles is performed by induction on $i$.
All steps of induction are similar to the first one.
Here we only show how to perform the first step and to construct the families
$\{v^0_\mu\times V^0_\mu\}_{\mu\in M_0}$ and
$\{u^1_\lambda\times U^1_\lambda\}_{\lambda\in\Lambda_1}$.

Put $M_0= X$ and for a point $x\in X$ consider a rectangle
$u^0_\lambda\times U_\lambda^0$ containing $\{x\}\times F(x)$.
By $UV^{[L]}$-property of $F(x)$ we find a closed neighborhood
$V^0_x$ of $F(x)$ such that the pair $V^0_x\subset U^0_\lambda$
is $[L]$-connected.
Since $F$ is u.s.c., we can choose a closed neighborhood
$v_x^0\subset u^0_\lambda$
of the point $x$ such that $F(v_x^0)\subset \Int V^0_x$.

Now we construct a family
$\{u^1_\lambda\times U^1_\lambda\}_{\lambda\in\Lambda_1}$.
Let $\alpha$ be a locally finite open cover of $X$ refining $v^0$.
For every element $A\in \alpha$ take an index $\mu\in M_0$ such that
$A\subset v^0_\mu$ and denote $W_A=\Int V^0_\mu$.
Then $A\times W_A$ lies in $v^0_\mu\times V^0_\mu$.
Let $u^1=\{u^1_\lambda\}_{\lambda\in\Lambda_1}$ be an open cover
of $X$ which is star-refined into $\alpha$. Define

\[  U^1_\lambda=\bigcap \{W_A\mid \St(u^1_\lambda,u^1)\subset A\in\alpha\}.  \]

To verify (5), consider $u^1_{\lambda'}\in u^1$ such that
$u^1_{\lambda'}\cap u^1_{\lambda}\ne\emptyset$.
Then $u^1_{\lambda'}\subset \St(u^1_\lambda,u^1) \subset A$ for some $A \in \alpha$ and by definition $U^1_{\lambda'}\subset W_A$.
Thus, $u^1_{\lambda'}\times U^1_{\lambda'}\subset A\times W_A
\subset v^0_\mu\times V^0_\mu$.
\bigskip

{\bf 2. Construction of the map $f$.}
\medskip

Since $X$ is a paracompact $C$-space, there exists a locally finite open cover
$\Sigma$ of $X$ of the form $\Sigma=\cup_{i=1}^\infty \sigma_i$ such that for
$i\ge 1$, $\sigma_i$ is pairwise disjoint collection refining $u^i$.
For every integer $k\ge 0$ denote by $\Sigma^{(k)}$ the set of points
$x\in X$ such that the cover $\Sigma$ has order $\le k+1$ at $x$.
Note that $X=\cup_{i=0}^\infty \Sigma^{(k)}$ and $\Sigma^{(k)}$
is closed in $X$. We will construct $f$ inductively
extending it over sets $\Sigma^{(k)}$.

For any element $s$ of the cover $\Sigma$ we denote by $\sigma(s)$
the integer number such that $s\in \sigma_{\sigma(s)}$.

{\it Choice 3:} for any element $s\in \Sigma$ we fix
$\lambda(s)\in \Lambda_{\sigma(s)}$
such that $s\subset u_{\lambda(s)}^{\sigma(s)}$.

Let $s_0, s_1,\dots,s_n$ be elements of the cover $\Sigma$
such that $s_0\cap s_1\cap\dots\cap s_n\ne\emptyset$.
Then this set of elements could be ordered according to
the order of numbers $\sigma(s_0),\sigma(s_1),\dots,\sigma(s_n)$,
and the smallest element of the set $\{s_0, s_1,\dots,s_n\}$
is called the {\it basic element}.
We always assume that $s_0$ is the basic element of
the set $\{s_0, s_1,\dots,s_n\}$.
We will use the following notations
\[ [s_0, s_1,\dots,s_n]=X \setminus \bigcup\{ \Sigma \setminus \{ s_{0},s_{1},\dots ,s_{n}\}\}  \]
\[ \langle s_0, s_1,\dots,s_n\rangle=(s_0\cap s_1\cap\dots\cap s_n)
   \cap \Sigma^{(n)}. \]
One should understand the set $[s_0, \dots,s_n]$ as closed
$n$-dimensional "simplex" with interior $\langle s_0, s_1,\dots,s_n\rangle$
and boundary $\cup_{m=0}^n [s_0,s_1,\dots,\wh s_m,\dots,s_k]$.
It is easy to check that
$\Sigma^{(n)}=\bigcup [s_{i_0},s_{i_1},\dots,s_{i_n}]$ and
\[ [s_0, \dots,s_n]= \cup_{m=0}^n [s_0,\dots,\wh s_m,\dots,s_k]
\cup \langle s_0, \dots,s_n\rangle.  \]

Let us construct the mapping $f$ on the set $\Sigma^{(0)}$
which is a discrete collection of sets of the type $[s_0]$.
We define $f$ independently on every such a set.
For a set $[s_0]$ we take a point $p\in F([s_0])$ and put $f([s_0])=p$.

Let us extend $f$ to arbitrary nonempty set $\langle s_0, s_1\rangle $.
For $i=0,1$ we have $\langle s_i\rangle \subset u^{\sigma(s_i)}_{\lambda(s_i)}$
and then $f(\langle s_i\rangle )\subset U^{\sigma(s_i)}_{\lambda(s_i)}$
by property (3). According to the choice 2, we take
$\mu\in M_{\sigma(s_0)-1}$ such that

$$  [s_0, s_1]\subset
    \overline{\St(u^{\sigma(s_0)}_{\lambda(s_0)},u^{\sigma(s_0)})}\subset
    v^{\sigma(s_0)-1}_\mu \qquad \text{and} \qquad
    f([s_0])\cup f([s_1])\subset V^{\sigma(s_0)-1}_\mu.   $$
Choice 1 gives us $\lambda=\lambda(\mu)$, a set $U^{\sigma(s_0)-1}_\lambda$
and $[L]$-extension

$$  e^{\sigma(s_0)-1}_\mu\colon \left( \wt V^{\sigma(s_0)-1}_\mu,
    W^{\sigma(s_0)-1}_\mu\right) \to \left(V^{\sigma(s_0)-1}_\mu,
    U^{\sigma(s_0)-1}_\lambda\right).   $$

Since the mapping $e^{\sigma(s_0)-1}_\mu|_{\wt V^{\sigma(s_0)-1}_\mu}$
is $[L]$-soft, we can lift the map
$f|_{[s_0]\cup [s_1]} \colon [s_0]\cup [s_1] \to V^{\sigma(s_0)-1}_\mu$
to a map $g\colon [s_0]\cup [s_1] \to \wt V^{\sigma(s_0)-1}_\mu$.
Now extend $g$ to a mapping
$\wt g\colon [s_0,s_1] \to W^{\sigma(s_0)-1}_\mu$
and define $f|_{[s_0,s_1]}$ as $e^{\sigma(s_0)-1}_\mu\circ \wt g$.

We can continue our construction so that the extension to a set
$\langle s_0, s_1,\dots,s_m\rangle $ uses $[L]$-extension $e^{\sigma(s_0)-1}_\mu$
and goes through $W^{\sigma(s_0)-1}_\mu$ resulting as
$f|_{[s_0,\dots,s_m]}=e^{\sigma(s_0)-1}_\mu\circ \wt g$.
Therefore, the set $f([s_0,\dots,s_m])$ is contained in
$U^{\sigma(s_0)-1}_\lambda$ while the set $[s_0,\dots,s_m]$
lies in $u^{\sigma(s_0)-1}_\lambda$.
Note that both indeces $\lambda$ and $\mu$ depend only on the basic
element $s_0$ and do not depend on $m$.
So, $[L]$-extension $e^{\sigma(s_0)-1}_\mu$ is a "rule"
for constructing mapping on each set $\langle s_0,\dots,s_m\rangle $
with basic element $s_0$.

Suppose that the map $f$ is constructed on $\Sigma^{(k-1)}$.
Let us extend $f$ independently to every set of type
$\langle s_0,\dots,s_k\rangle $.
Since the difference $\Sigma^{(k)} \setminus \Sigma^{(k-1)}$
is covered by a discrete family of such sets, it follows that
the so obtained extension of $f$ to $\Sigma^{(k)}$ would be continuous.
Assume that $s_1$ is basic element of the set $\{s_1, s_2,\dots,s_k\}$.
Then the set $f(\langle s_1,\dots,s_k\rangle )$ lies in some $U^{\sigma(s_1)-1}_{\lambda_1}$
and $u^{\sigma(s_1)-1}_{\lambda_1}$ contains $[s_1,\dots,s_k]$.
Since $\sigma(s_1)-1\ge \sigma(s_0)$, the set $f(\langle s_1,\dots,s_k\rangle )$
lies in $V^{\sigma(s_0)-1}_{\mu}$ by property (5). Let

\[  G=\bigcup_{1\le m\le k} [s_0,s_1,\dots,\wh s_m,\dots,s_k]. \]

\noindent Then, by our construction, $f|_G$ has a lift
$g\colon G\to W^{\sigma(s_0)-1}_{\mu}$.
Note that

\[ f(\overline{\langle s_1,\dots,s_k\rangle }\cap G) \subseteq V^{\sigma(s_0)-1}_{\mu} = \overline{V^{\sigma(s_0)-1}_{\mu}} .\]

\noindent Since the mapping
$e^{\sigma(s_0)-1}_\mu\colon \wt V^{\sigma(s_0)-1}_\mu
\to V^{\sigma(s_0)-1}_\mu$ is $[L]$-soft, we extend the lift $g$
to the set $\overline{\langle s_1,\dots,s_k\rangle}$.
Now extend it to a mapping
$g\colon [s_0,\dots,s_k] \to W^{\sigma(s_0)-1}_\mu$
and define $\displaystyle f|_{[s_0,\dots,s_k]}$ as the
composition $\displaystyle e^{\sigma(s_0)-1}_\mu\circ g$.

It only remains to note that the local finiteness of $\Sigma$
guarantees the continuity of the above constructed map $f$.
Proof is completed.

In conclusion authors would like to thank the referee for several helpful suggestions.

\end{document}